\documentclass[12pt]{article}
\usepackage{latexsym}
\usepackage{amsmath}
\usepackage{amssymb}
\usepackage{amsfonts}
\usepackage{graphics,graphicx}

\usepackage{url}

\def\claim#1{\begin{trivlist}\item[\hskip\labelsep\bf#1]\it}
\def\endclaim{\end{trivlist}}

\numberwithin{equation}{section}

\newtheorem{theorem}{Theorem}[section]
\newtheorem{lemma}[theorem]{Lemma}

\newcommand{\pf}{\noindent{\bf Proof}}

\title{Unique continuation property for multi-terms time fractional diffusion equations}
\author{Ching-Lung Lin\thanks{Department of Mathematics, NCTS, National Cheng-
Kung University, Tainan 701, Taiwan. Partially supported by the
Ministry of Science and Technology of Taiwan. (Email:
cllin2@mail.ncku.edu.tw)}\qquad Gen Nakamura\thanks{Department of
Mathematics, Hokkaido University, Sapporo 060-0808, Japan. Partially supported by
grant-in-aid for Scientific Research (15K21766 and 15H05740) of the Japan Society for the Promotion of Science. (Email: gnaka@math.sci.hokudai.ac.jp)
}}

\date{}

\begin{document}
\renewcommand{\theequation}{\thesection.\arabic{equation}}
 \maketitle
\begin{abstract}
A Carleman estimate and the unique continuation property of solutions for a multi-terms time fractional diffusion equation up to order $\alpha\,\,(0<\alpha<2)$
and general time dependent second order strongly elliptic time elliptic operator for the diffusion.
The estimate is derived via some subelliptic estimate for an operator associated to this equation using calculus of pseudo-differential operators.
A special Holmgren type transformation which is linear with respect to time is used to show the local unique continuation of solutions. 
We developed a new argument to derive the global unique continuation of solutions. Here the global unique continuation means as follows. 
If u is a solution of the multi-terms time fractional diffusion equation in a domain over the time interval $(0,T)$, 
then a zero set of solution over a subdomain of $\Omega$ can be continued to $(0,T)\times\Omega$.

\end{abstract}
\section{Introduction}\label{sec1}
Let $(0,T)$ be a time interval and $\Omega$ be a domain in $n$ dimensional Euclidean space ${\mathbf{R}}^n$.
Let $m\in {\mathbb N}$ and $\alpha_j, q_j$ for $j=1,2\cdots m$ be positive constants such that $q_1=1$ and $2>\alpha=\alpha_1>\alpha_2>\cdots>\alpha_m>0$.
Consider an anomalous diffusion equation with multi-terms fractional time derivatives given by
\begin{equation}\label{anomalous eq}
\Sigma_{j=1}^mq_j\partial_t^{\alpha_j} u(t,y)-\emph{L} u(t,y)=\ell_1(t,y;\nabla_y)u(t,y),\,\,t\in(0,T),\,y\in\Omega,
\end{equation}
where the fractional derivative $\partial_t^\alpha u$ of order $0<\alpha<1$ or $1<\alpha<2$ in the Caputo sense is given by
\begin{equation}\label{1.2}
\begin{array}{l}
\partial_t^\alpha
=\frac{1}{\Gamma(k-\alpha)}\big((t^{k-1-\alpha}H(t)\otimes\delta_y)*\partial^k_t u\big)(t,y)\quad k-1<\alpha<k, k=1,2,
\end{array}
\end{equation}
with the Heaviside function $H(t)$, the Dirac delta function $\delta_y$ supported at $y$ and
a linear partial differential operator $\ell_1(t,y;\nabla_y)$ of order $1$.
We assume the coefficients of $\ell_1(t,y;\nabla_y)$ are $C^\infty([0,T]\times\overline\Omega)$.
The coefficients could be just in $C^0([0,T]\times\overline\Omega)$, but without loss of generality
and simplicity we have assumed in this way. Also $\emph{L}$ is defined by
$$\emph{L} u(t,y):=\sum_{j,k=1}^na_{jk}(t,y)\partial_{y_j}\partial_{y_k}u(t,y).$$
Here we assume that $a_{jk}(t,y)=a_{kj}(t,y)\in C^\infty([0,T]\times\overline\Omega)$
and the estimate
\begin{equation}\label{1.3}
\begin{array}{l}
\delta |\xi|^2 \leq \sum_{j,k=1}^na_{jk}(t,y)\xi_j\xi_k \leq \delta^{-1} |\xi|^2, \quad \xi=(\xi_1,\cdots,\xi_n)
\end{array}
\end{equation}
holds for some constant $\delta>0$.
In the sequel we always extend each $a_{jk}\in C^\infty([0,T]\times\overline\Omega)$ to
$a_{jk}\in B^\infty(\mathbf{R}^{1+n})$ without destroying \eqref{1.3} and $a_{jk}(t,y)=a_{kj}(t,y)$, where
$B^\infty(\mathbf{R}^{1+n})$ is the set of all functions in $C^\infty(\mathbf{R}^{1+n})$ which
are bounded together with their derivatives.

\medskip

The aim of this paper is to give the unique continuation property abbreviated by UCP for solutions of \eqref{anomalous eq}. More precisely we have the following main theorem.

\begin{theorem}\label{thm1.1}{\rm(global UCP)}
For any solution $u\in H^{\alpha,2}((0,T)\times\Omega)$
(see below in this section for its definition) of \eqref{anomalous eq}
supported on $\{t\ge0\}$, if $u=0$ in a subdomain of $\Omega$ over $(0,T)$, then
$u=0$ in $(0,T)\times\Omega$.
\end{theorem}

This will follows from the following theorem.

\begin{theorem}\label{thm1.2}{\rm (local UCP)}
Let $0<T'<T$, $\hat y=(\hat{y}_1,\cdots,\hat{y}_{n-1},0)\in\Omega$,
$\omega=\{(y_1,\cdots,y_n):\,\hat{y}_j-\ell<y_j<\hat{y}_j+\ell,\,1\le j\le n-1,\,-l<y_n< 0\}\subset\Omega$ with $\ell>0$.
Assume that  $u(t,y)\in H^{\alpha,2}((0,T)\times\omega)$ satisfies
\begin{equation}\label{1.4}
\left\{
\begin{array}{l}
\Sigma_{j=1}^mq_j\partial_t^{\alpha_j} u(t,y)-\emph{L} u(t,y)=\ell_1(t,y;\nabla_y)u(t,y),\\
u(t,y)=0,\,\,t\le0,\\
u(t,y)=0,\,\,y\in\omega,\,t\in(0,T'),
\end{array}
\right.
\end{equation}
Then $u$ will be zero in a neighborhood $V(T',\omega)$ of $(t,y)=(0,\hat y)$.
\end{theorem}

Concerning the physical background of the anomalous diffusion equation, this equation with single-term fractional derivative was first studied in material science and the exponent $\alpha$ of $\partial^{\alpha}_{t}u$ in this equation is an index which describes the long time behavior of the mean square displacement $<x^2(t)>\sim\mbox{\it positive const.}\,t^\alpha$ of a diffusive particle $x(t)$ describing anomalous diffusion on fractals such as some amorphous semiconductors or strongly
porous materials (see \cite{Anh}, \cite{Metzler} and references therein).

Recently there are more and more studies which use anomalous diffusion equations as model equations to describe physical pheonmena. For instance in enviormental science, an anomalous diffusion equation was used to describe the spread of pollution in soil which cannot correctly modeled by the ususal diffusion equationa (see \cite{hat1}, \cite{hat2}).

The Cauchy problem and intial boundary value problem for the anomalous diffusion equation have been already studied by many people (see \cite{A}, \cite{E} and the references therein). Then, this recent research trend triggered some mathemtatical studies of inverse problems on anomalous diffusion equations (see \cite{Jin_Rundell}, \cite{Li_Yamamoto}, \cite{Miller_Yamamoto}).

One very important key to the study of control problem and inverse problem is UCP for a model equation of the problem. It can give the approximate boundary controllability for the control problem and for reconstruction schemes such as the linear sampling method to identify unknown objects such as cracks, cavities and inclusions inside an anomalous diffusive medium. The usual way to have UCP is by a Carleman estimate. As for the time fractional diffusion equations, Carleman estimates have been given for some special cases. That is for the single-term time fractional diffusion equation,
$\alpha=1/2$, a Carleman estimate was given in \cite{XCY}, \cite{ZX}, \cite{YZ} for $n=1$ and \cite{CLN} for $n=2$ via that for the operator $\partial_t+\Delta^2$ with some
lower order terms. Further for the single-term time fractional diffusion equation with second order elliptic operator with homogeneous isotropic diffusion coefficient, a local Carleman estimate giving the global UCP was given in \cite{LN} when the order $\alpha$ of time fractional derivative is $0<\alpha<1$.

We can basically say that the the multi-terms time fractional diffusion equation is semi-elliptic. Treve (\cite{Treves}) gave a way to derive a Carleman estimate which can have UCP for some semi-elliptic equations. We will adapt Treves' argument to have UCP for the multi-term time fractional diffusion equation conjugated by $e^{t}$ and transformed by a Holmgren transform which is still subelliptic. Our method is quite general and has a potential to be applied having UCP to space time fractional diffusion equations (\cite{Hanyga}) and some fractional derivative visco-elastic equations (\cite{Lu}).

The rest of this paper is organized as follows.
In Section 2, we compute the principal part of the anomalous diffusion operator conjugated by $e^t$
which undergone Holmgren type transformation. Section 3 is devoted to
analyzing the Poisson bracket of its principal symbol.
In Section 4, we derive some subelliptic estimate for some pseudo-differential operator
associated to this principal part of the operator. Then by using this subelliptic estimate,
the Carleman estimate is derived in Section 5. In the remaining two sections we give the UCP.
More precisely in Section 6, the local UCP stated in Theorem \ref{thm1.2} is proved, and in Section 7,
the global UCP stated in Theorem \ref{thm1.1} is proved.

We close this section by giving the definition of Sobolev space
$H^{\alpha,2}({\mathbf{R}}\times\Omega)$.
Let $\mathcal{S}(\mathbf{R}_t^1\times\mathbf{R}_x^n)$ and $\mathcal{S}'(\mathbf{R}_t^1\times\mathbf{R}_x^n)$ be
the set of rapidly decreasing functions in $\mathbf{R}_t^1\times\mathbf{R}_x^n$ and its dual space, respectively.
Then, for $m,s\in \mathbf{R}$, $v=v(t,x)\in \mathcal{S}'(\mathbf{R}_t^1\times\mathbf{R}_x^n)$, belongs to the function space
$H^{m,s}(\mathbf{R}_t\times\mathbf{R}_x^n)$ if
$$\|v\|^2_{H^{m,s}}:=\iint(1+|\xi|^s+|\tau|^m)^2|\hat{v}|^2d\tau d\xi$$
is finite, and $\|v\|^2_{H^{m,s}}$ denotes the norm of $v$ of this function space, where
$\hat{v}$ is the Fourier transform defined by
$$\mathcal{F}(v)(\tau,\xi)=\hat{v}(\tau,\xi)=\iint e^{-it\tau-ix\cdot\xi}v(t,x)dxdt.$$
Further, for any open sets $A\subset{\mathbf{R}}_t,\,B\subset{\mathbf{R}}_x^n$,
we define $H^{m,s}(A\times B)$ as the restriction of $H^{m,s}(\mathbf{R}_t\times\mathbf{R}_x^n)$ to $A\times B$.

\section{Holmgren type transformation and $\Lambda_\alpha^m(D_t,D_x)$}\label{sec2}
\setcounter{equation}{0}
In this section, we prepare some setups which are necessary to prove Theorem \ref{thm1.2}.
In particular, they are Holmgren type transformation and $\Lambda_\alpha^m(D_t,D_x)$.

Using a cutoff function on $y$ and a coordinate transformation in $y$ coordinates, we can consider our problem in $|y-\hat{y}|$ small for the aforementioned $\hat{y}\in \mathbf{R}^n$ and assume that
$$u(t,y)=0 \quad (y_n<0\quad {\rm or}\quad t<0).$$
To be precise, we define for $\zeta\in\mathbf{N}$ that
\begin{equation*}
\begin{array}{l}
\begin{cases}
Q_\zeta=\{(y_1,\cdots,y_n):|y_j-\hat{y}_j|\leq \sqrt{X},-l/3<y_n\leq \zeta X,1\leq j\leq n-1 \},\\
\tilde{Q}_\zeta=\{(y_1,\cdots,y_n):|y_j-\hat{y}_j|\leq 2\sqrt{X},-2l/3<y_n\leq (\zeta+1) X, 1\leq j\leq n-1\},
\end{cases}
\end{array}
\end{equation*}
where $X$ is a small positive constant which will be determined later. Let us define
$\kappa_\zeta(y)\in C^\infty(\mathbf{R}^n)$ by
\begin{equation*}
\begin{array}{l}
\kappa_\zeta(y)=
\begin{cases}
1, \quad y\in Q_\zeta,\\
0, \quad y\in \mathbf{R}^n\setminus \tilde{Q}_\zeta.
\end{cases}
\end{array}
\end{equation*}
We clearly have
$$(\kappa_1 u)(t,y)=0,\quad \partial^{\alpha}_{t}(\kappa_1 u)(t,y)=0 \quad (y_n \leq 0\quad {\rm or}\quad t \leq 0 \quad {\rm or}\quad y\in \mathbf{R}^n\setminus \tilde{Q}_1)$$
and
$$(\kappa_1 u)(t,y)=u(t,y), \quad y\in Q_1.$$

To prove Theorem \ref{thm1.2},
we use the change of variables
\begin{equation}\label{1.45}
\begin{array}{l}
x'=y'-\hat{y}', x_n=y_n+c|y'-\hat{y}'|^2+\frac{X}{T}t, \tilde{t}=t,
\end{array}
\end{equation}
where $x'=(x_1,\cdots,x_{n-1})$, $y'=(y_1,\cdots,y_{n-1})$, $\hat{y}'=(\hat{y}_1,\cdots,\hat{y}_{n-1})$, $c\geq 1$ is a suitable positive constant and $X$ is a small positive constant which will be determined later.
This is a Holmgren type transformation. We sometimes called $(t,x)$ the Holmgren coordinates.
By $\partial_t=\frac{X}{T}\partial_{x_n}+\partial_{\tilde{t}}$,
$ \partial_{y_j}=\partial_{x_j}+2cx_j\partial_{x_n}$ for $j=1,\cdots,n-1$ and $\partial_{y_n}=\partial_{x_n}$, we have
\begin{equation}\label{1.5}
\begin{array}{l}
\partial_{t}^\alpha u(t,y)=\frac{1}{\Gamma(1-\alpha)}\int_0^t(t-\eta)^{-\alpha}\partial_\eta u(\eta,y)d\eta\\
=\frac{1}{\Gamma(1-\alpha)}\int_0^{\tilde{t}}(\tilde{t}-\tilde{\eta})^{-\alpha}\partial_{\tilde{\eta}} u(\tilde{\eta},x)d\tilde{\eta}
+\frac{X}{T\Gamma(1-\alpha)}\int_0^{\tilde{t}}(\tilde{t}-\tilde{\eta})^{-\alpha}\partial_{x_n}u(\tilde{\eta},x)d\tilde{\eta},
\end{array}
\end{equation}
where $\Gamma(\cdot)$ denotes the Gamma function and $u(\tilde{t},x)$
with $\tilde{t}=\tilde{\eta}$ on the right hand side of \eqref{2.2}
is the push forward of $u(t,y)$ by the change of variables.
Since
$$(\kappa_1 u)(t,y)=0,\quad \partial^{\alpha}_{t}(\kappa_1 u)(t,y)=0 \quad (y_n \leq 0\quad {\rm or}\quad t \leq 0 \quad {\rm or}\quad y\in \mathbf{R}^n\setminus \tilde{Q}_1)$$
and
$x_n=y_n+c|y'-\hat{y}'|^2+\frac{X}{T}t$, $u(\tilde{t},x)$ satisfies
$${\rm supp}\, \kappa_1 u\subset\{x_n\geq 0\}.$$
Further, let $\chi$ be a smooth function such that
\begin{equation}\label{1.6}
\chi (x_n)=
\begin{cases}
\begin{array}{l}
1,\quad x_n\leq (1-\epsilon)X,\\
0,\quad x_n\geq X.
\end{array}
\end{cases}
\end{equation}
Then
\begin{equation}\label{1.7}
{\rm supp}\, \chi\kappa_1 u\subset Q_1.
\end{equation}
Since $\kappa_1 u=u$ in $Q_1$, we will only denote $\kappa_1 u$ by $u$. Also, we will assume that
\begin{equation}\label{1.8}
|x'|\lesssim \sqrt{X}, \quad x_n\leq X,
\end{equation}
here and after $A \lesssim B$ denotes $A \leq CB$ for positive constant $C$ depending on $n$ and $\alpha$.

For further arguments, we define a function space $\dot{H}^m_{\alpha}(\overline{\mathbf{R}_+^{1+n}})$ for $m\in \mathbf{R}$ as follows.
First, let $(t,x)\in\mathbf{R}^{1+m}$ and denote the open half of $\mathbf{R}^{1+m}$ space defined by $t>0$ and complement of its closure $\bar{\mathbf{R}}_+^{1+m}$
by $\mathbf{R}_+^{1+m}$ and $\mathbf{R}_-^{1+m}$, respectively. If $E$ is a space of distributions in $\mathbf{R}^{1+m}$ we
use the notation $\bar{E}(\mathbf{R}_+^{1+m})$ for the space of restrictions to $\mathbf{R}_+^{1+m}$ of elements in $E$ and we write
$\dot{E}(\bar{\mathbf{R}}_+^{1+m})$ for the set of distributions in $E$ supported by $\bar{\mathbf{R}}_+^{1+m}$.

Our strategy is to use the symbol calculus. To begin with set
$$\Lambda_{\alpha}^m(\tau,\xi)=((1+|\xi|^2)^{1/\alpha}+i\tau)^{m\alpha/2} \quad \mbox{for}\,\,m\in \mathbf{R}$$
and define a pseudo-differential operator $\Lambda_{\alpha}^m(D_t,D_x)$ by $\Lambda_{\alpha}^m(D_t,D_x)\phi:=\mathcal{F}^{-1}(\Lambda_{\alpha}^m(\tau,\xi)\hat{\phi})$, $\phi\in \dot{\mathcal{S}}(\overline{\mathbf{R}_+^{1+n}})$, where
$\mathcal{F}^{-1}$ is the inverse Fourier transform, $S$, $S'$ are space of the Schwartz space and
its dual space in $\mathbf{R}^{1+n}$. Then
$\Lambda_{\alpha}^m(D_t,D_x)$ is a continuous map on $\dot{\mathcal{S}}(\overline{\mathbf{R}_+^{1+n}})$
with inverse $\Lambda_{\alpha}^{-m}(D_t,D_x)$. (See Theorem B.2.4 in \cite{Hormander}.) Based on this we define $\dot{H}^m_{\alpha}(\overline{\mathbf{R}_+^{1+n}})$ by $\dot{H}^m_{\alpha}(\overline{\mathbf{R}_+^{1+n}}):=\Lambda_{\alpha}^{-m}L^2(\mathbf{R}_+^{1+n})$,
where we considered $L^2(\mathbf{R}_+^{1+n})$ as a subset of $\dot{\mathcal{S}}'(\overline{\mathbf{R}_+^{1+n}})$.

To avoid having zero for the fractional power of $\tau$, we conjugate $\sum_{l=1}^m\partial_t^{\alpha_l}-\emph{L}$ with
$e^{t}$. Thus let $P(t,x,D_t,D_x)$ be the operator $e^{-t}(\sum_{l=1}^m\partial_t^{\alpha_l}-\emph{L})e^{t}$
in terms of the Holmgren coordinates $(t,x)$ (see \eqref{1.45}) and
$p(t,x,\tau,\xi)$ be its total symbol, where $D_t=-\sqrt{-1}\partial_t$, $D_x=-\sqrt{-1}\partial_x$. Then it is easy to see
\begin{equation}\label{2.1}
\begin{array}{rl}
p(t,x,\tau,\xi)=&\Sigma_{l=1}^mq_l(1+i\tau)^{\alpha_l}+a_{nn}(t,x)\xi_n^2+2\Sigma_{j=1}^{n-1}a_{jn}(t,x)\xi_n(\xi_j+2cx_j\xi_n)\\
&+\Sigma_{j,k=1}^{n-1}a_{jk}(t,x)(\xi_j+2cx_j\xi_n)(\xi_k+2cx_k\xi_n)\\
&+\Sigma_{l=1}^mq_l\frac{X}{T}i^\alpha(\tau-i)^{\alpha_l-1}\xi_n.
\end{array}
\end{equation}
We want to derive a Carleman estimate for the operator $P(t,x,D_t,D_x)$
using some subelliptic estimate based on the idea by F. Treves (\cite{Treves}).

\section{Poisson bracket}\label{sec3}
\setcounter{equation}{0}
First, express $(1+i\tau)^\alpha$ in the form
$$
(1+i\tau)^\alpha=|1+i\tau|^\alpha[\cos(\alpha \arg(1+i\tau))+i\sin(\alpha \arg(1+i\tau))].
$$
The crucial observation is that there exists a positive constant $C_\alpha=\min\{\frac{\sqrt{2}}{2},\sin(\pi(1-\frac{\alpha}{2}))\}$ such that
\begin{equation}\label{2.2}
|\Im (1+i\tau)^\alpha|\geq C_\alpha|1+i\tau|^\alpha,\qquad |\tau|\geq 1.
\end{equation}
We also have
\begin{equation*}
|(1+i\tau)^{\alpha-1}|\leq|1+i\tau|^{\alpha-1}.
\end{equation*}
Based on this, in terms of the scaling
$$
(\xi,\tau)\mapsto (\rho\xi,\rho^{2/\alpha}\tau)
$$
with large $\rho>0$, the principal part $\tilde{p}(t,x,\tau,\xi)$ of  $p(t,x,\tau,\xi)$ is
\begin{equation*}
\begin{array}{rl}
\tilde{p}(t,x,\tau,\xi)=&(1+i\tau)^{\alpha}+a_{nn}(t,x)\xi_n^2+2\Sigma_{j=1}^{n-1}a_{jn}(t,x)\xi_n(\xi_j+2cx_j\xi_n)\\
&+\Sigma_{j,k=1}^{n-1}a_{jk}(t,x)(\xi_j+2cx_j\xi_n)(\xi_k+2cx_k\xi_n)\\
\end{array}
\end{equation*}
Let $\psi=\frac{1}{2}(x_n-2X)^2$ and consider the symbol $p(t,x,\tau,\xi+i|\sigma|\nabla\psi)$ over $\mathbf{R}^{n+1}\times\mathbf{R}_z$ to define a pseudo-differential operator $P_\psi(t,x,D_t,D_{x},D_z)$ by
$$P_\psi=P_\psi(t,x,D_t,D_{x},D_z)=p(t,x,D_{t},D_x+i|D_z|\nabla \psi)$$
which is given for any compactly supported distribution $v$ in $\mathbf{R}^{n+1}\times\mathbf{R}_z$ by
\begin{equation}\label{2.3}
\begin{array}{rl}
P_\psi v(z,t,x)=\int e^{i(x\cdot\xi+t\tau+z\sigma)}p(t,x,\tau,\xi+i|\sigma|\nabla\psi)\hat{v}(\sigma,\xi,\tau)d\sigma d\tau d\xi.
\end{array}
\end{equation}
We denote the principal symbol of $P_\psi$ by $\tilde{p}_\psi$ which is given by
\begin{equation*}
\begin{array}{rl}
\tilde{p}_\psi=&(1+i\tau)^{\alpha}+a_{nn}(\xi_n+i|\sigma|\tilde{X})^2\\
&+2\Sigma_{j=1}^{n-1}a_{jn}(\xi_j+2cx_j\xi_n+2icx_j|\sigma|\tilde{X})(\xi_n+i|\sigma|\tilde{X})\\
&+\Sigma_{j,k=1}^{n-1}a_{jk}(\xi_j+2cx_j\xi_n+2icx_j|\sigma|\tilde{X})(\xi_k+2cx_k\xi_n+2icx_k|\sigma|\tilde{X}),
\end{array}
\end{equation*}
where $\tilde{X}=x_n-2X$.

By the definition of Poisson bracket
\begin{equation}\label{2.4}
\{\Re \tilde{p}_\psi, \Im \tilde{p}_\psi\}=
\Sigma_{j=1}^{n}(\partial_{\xi_j}\Re \tilde{p}_\psi\cdot\partial_{x_j}\Im \tilde{p}_\psi-
\partial_{x_j}\Re \tilde{p}_\psi\cdot\partial_{\xi_j}\Im \tilde{p}_\psi)
+\partial_{\tau}\Re \tilde{p}_\psi\cdot\partial_{t}\Im \tilde{p}_\psi-
\partial_{t}\Re \tilde{p}_\psi\cdot\partial_{\tau}\Im \tilde{p}_\psi
\end{equation}
with the real part $\Re \tilde{p}_\psi$ and imaginary part $\Im \tilde{p}_\psi$ of $\tilde{p}_\psi$.
The principal part $\{\Re \tilde{p}_\psi, \Im \tilde{p}_\psi\}_p$ of the Poisson bracket
$\{\Re \tilde{p}_\psi, \Im \tilde{p}_\psi\}$ is
\begin{equation*}
\begin{array}{ll}
\{\Re \tilde{p}_\psi, \Im \tilde{p}_\psi\}_p=
\Sigma_{j=1}^{n}(\partial_{\xi_j}\Re \tilde{p}_\psi\cdot\partial_{x_j}\Im \tilde{p}_\psi
-\partial_{x_j}\Re \tilde{p}_\psi\cdot\partial_{\xi_j}\Im \tilde{p}_\psi).
\end{array}
\end{equation*}
Write $\tilde{p}_\psi=\tilde{p}_{1,\psi}+\tilde{p}_{2,\psi}$, where
\begin{equation*}
\begin{array}{rl}
\tilde{p}_{1,\psi}=&a_{nn}(\xi^2_n-|\sigma|^2\tilde{X}^2)+2\Sigma_{j=1}^{n-1}a_{jn}\xi_j\xi_n+\Sigma_{j,k=1}^{n-1}a_{jk}\xi_j\xi_k\\
&+i[2a_{nn}\xi_n|\sigma|\tilde{X}+2\Sigma_{j=1}^{n-1}a_{jn}\xi_j|\sigma|\tilde{X}]
\end{array}
\end{equation*}
and
\begin{equation*}
\begin{array}{rl}
\tilde{p}_{2,\psi}=&(1+i\tau)^{\alpha}+2\Sigma_{j=1}^{n-1}a_{jn}(2cx_j\xi_n+2icx_j|\sigma|\tilde{X})(\xi_n+i|\sigma|\tilde{X})\\
&+\Sigma_{j,k=1}^{n-1}a_{jk}(2cx_j\xi_n+2icx_j|\sigma|\tilde{X})(\xi_k+2cx_k\xi_n+2icx_k|\sigma|\tilde{X})\\
&+2\Sigma_{j,k=1}^{n-1}a_{jk}\xi_j(2cx_k\xi_n+2icx_k|\sigma|\tilde{X}).
\end{array}
\end{equation*}

A direct computation gives
\begin{equation*}
\begin{array}{rl}
\partial_{\xi_n}\Re \tilde{p}_{1,\psi}=&2\Sigma_{j=1}^{n}a_{jn}\xi_j\\
\partial_{x_n}\Im \tilde{p}_{1,\psi}=&2\partial_{x_n}(\Sigma_{j=1}^{n}a_{jn}\xi_j)|\sigma|\tilde{X}+2|\sigma|\Sigma_{j=1}^{n}a_{jn}\xi_j\\
=&2|\sigma|\Sigma_{j=1}^{n}a_{jn}\xi_j+O(|\xi|^2+|\sigma|^2X^2)\\
\partial_{x_n}\Re \tilde{p}_{1,\psi}=&(\xi^2_n-|\sigma|^2\tilde{X}^2)\partial_{x_n}(a_{nn})+2\Sigma_{j=1}^{n-1}\partial_{x_n}(a_{jn})\xi_j\xi_n\\
&+\Sigma_{j,k=1}^{n-1}\partial_{x_n}(a_{jn})\xi_j\xi_k-2a_{nn}|\sigma|^2\tilde{X}\\
=&-2a_{nn}\sigma^2\tilde{X}+O(|\xi|^2+|\sigma|^2X^2)\\
\partial_{\xi_n}\Im \tilde{p}_{1,\psi}=&2a_{nn}|\sigma|\tilde{X}
\end{array}
\end{equation*}
and
\begin{equation*}
\begin{array}{rl}
\nabla_{\xi'}\Re \tilde{p}_{1,\psi}=&O(|\xi|+|\sigma|X)\\
\nabla_{x'}\Im \tilde{p}_{1,\psi}=&O(|\xi|^2+|\sigma|^2X^2)\\
\nabla_{x'}\Re \tilde{p}_{1,\psi}=&O(|\xi|^2+|\sigma|^2X^2)\\
\nabla_{\xi'}\Im \tilde{p}_{1,\psi}=&O(|\xi|+|\sigma|X).
\end{array}
\end{equation*}
Since $\tilde{p}_{2,\psi}$ contains $x_j$ in each terms, the derivation on $a_{jk}$ will preserve the order.
Thus, we have by \eqref{1.8}
\begin{equation}\label{2.5}
\begin{array}{rl}
\partial_{\xi_n}\Re \tilde{p}_{\psi}=&2\Sigma_{j=1}^{n}a_{jn}\xi_j+O(\sqrt{X}(|\xi|+|\sigma|X))\\
\partial_{x_n}\Im \tilde{p}_{\psi}=&2|\sigma|\Sigma_{j=1}^{n}a_{jn}\xi_j+O(\frac{1}{\sqrt{X}}(|\xi|+|\sigma|X)^2)\\
\partial_{x_n}\Re \tilde{p}_{\psi}=&-2a_{nn}\sigma^2\tilde{X}+O(\frac{1}{\sqrt{X}}(|\xi|+|\sigma|X)^2)\\
\partial_{\xi_n}\Im \tilde{p}_{\psi}=&2a_{nn}|\sigma|\tilde{X}+O(\sqrt{X}(|\sigma|X))
\end{array}
\end{equation}
and
\begin{equation}\label{2.6}
\begin{array}{rl}
\nabla_{\xi'}\Re \tilde{p}_{\psi}=&O(|\xi|+|\sigma|X)\\
\nabla_{x'}\Im \tilde{p}_{\psi}=&O((|\xi|+|\sigma|X)^2)\\
\nabla_{x'}\Re \tilde{p}_{\psi}=&O((|\xi|+|\sigma|X)^2)\\
\nabla_{\xi'}\Im \tilde{p}_{\psi}=&O(|\xi|+|\sigma|X).
\end{array}
\end{equation}
Combining \eqref{2.5} and \eqref{2.6}, we have
\begin{equation}\label{2.7}
\begin{array}{rl}
\{\Re \tilde{p}_\psi, \Im \tilde{p}_\psi\}_p=&\Sigma_{j=1}^{n}(\partial_{\xi_j}\Re \tilde{p}\cdot\partial_{x_j}\Im \tilde{p}
-\partial_{x_j}\Re \tilde{p}\cdot\partial_{\xi_j}\Im \tilde{p})\\
=&4|\sigma|(\Sigma_{j=1}^{n}a_{jn}\xi_j)^2+4a_{nn}^2|\sigma|^3\tilde{X}^2+O(\frac{1}{\sqrt{X}}(|\xi|+|\sigma|X)^3)
\end{array}
\end{equation}
\begin{lemma}\label{lem2.1}
Let $\tilde{p}_{\psi}=0$ and $x$ satisfies \eqref{1.7}, then we have
\begin{equation}\label{2.8}
\begin{array}{l}
\{\Re \tilde{p}_\psi, \Im \tilde{p}_\psi\}_p \gtrsim (|\xi|^2+\sigma^2+|\tau|^{\alpha})^{3/2}.
\end{array}
\end{equation}
\end{lemma}
\pf. Recall that
\begin{equation*}
\begin{array}{rl}
\tilde{p}_\psi=&\Sigma_{l=1}^mq_l(1+i\tau)^{\alpha_l}+a_{nn}(\xi_n+i|\sigma|\tilde{X})^2\\
&+2\Sigma_{j=1}^{n-1}a_{jn}(\xi_j+2x_j\xi_n+2ix_j|\sigma|\tilde{X})(\xi_n+i|\sigma|\tilde{X})\\
&+\Sigma_{j,k=1}^{n-1}a_{jk}(\xi_j+2x_j\xi_n+2ix_j|\sigma|\tilde{X})(\xi_k+2x_k\xi_n+2ix_k|\sigma|\tilde{X})\\
=&\Sigma_{l=1}^mq_l(1+i\tau)^{\alpha_l}+\Sigma_{j,k=1}^{n}a_{jk}\xi_j \xi_k - a_{nn}|\sigma|^2\tilde{X}^2\\
&+2i|\sigma|\tilde{X}\Sigma_{j=1}^{n}a_{jn}\xi_j+O(\sqrt{X}(|\xi|+|\sigma|X)^2),
\end{array}
\end{equation*}
where $\tilde{X}=x_n-2X$.

Let us divide $\tau$ into two cases. If $|\tau| \leq 1$, then
\begin{equation}\label{2.9}
\begin{array}{l}
\Re (1+i\tau)^\alpha \geq \varepsilon_0 |1+i\tau|^\alpha
\end{array}
\end{equation}
where $\varepsilon_0$ is a positive constant.
From $\Re \tilde{p}_\psi=0$ and \eqref{2.9}, we have
\begin{equation*}
\begin{array}{l}
\Re\big(\Sigma_{l=1}^mq_l(1+i\tau)^{\alpha_l}\big)+\Sigma_{j,k=1}^{n}a_{jk}\xi_j \xi_k=a_{nn}|\sigma|^2\tilde{X}^2+O(\sqrt{X}(|\xi|+|\sigma|X)^2)
\end{array}
\end{equation*}
which implies
\begin{equation*}
\begin{array}{l}
|\xi|^2+|1+i\tau|^{\alpha}\lesssim |\sigma|^{2}X^2.
\end{array}
\end{equation*}

If  $|\tau| \geq 1$, then
\begin{equation}\label{2.10}
\begin{array}{l}
|\Im (1+i\tau)^{\alpha}| \geq \varepsilon_1 |1+i\tau|^{\alpha}
\end{array}
\end{equation}
where $\varepsilon_1$ is a positive constant.

First, we let
$$(2\Sigma_{l=1}^mq_l)|1+i\tau|^{\alpha} \leq \delta_0|\xi|^2,$$
where $\delta_0$ is a positive constant satisfying $\delta_0|\xi|^2\leq\Sigma_{j,k=1}^{n}a_{jk}\xi_j \xi_k$.
From  $\Re \tilde{p}_\psi=0$, we have
\begin{equation*}
\begin{array}{l}
\Re\big(\Sigma_{l=1}^mq_l(1+i\tau)^{\alpha_l}\big)+\Sigma_{j,k=1}^{n}a_{jk}\xi_j \xi_k=a_{nn}|\sigma|^2\tilde{X}^2+O(\sqrt{X}(|\xi|+|\sigma|X)^2)
\end{array}
\end{equation*}
which implies
\begin{equation}\label{2.11}
\begin{array}{l}
|\xi|^2+|(1+i\tau)|^{\alpha}\lesssim |\sigma|^{2}X^2.
\end{array}
\end{equation}

On the other hand, we let
$$\delta_0|\xi|^2 \leq (2\Sigma_{l=1}^mq_l)|1+i\tau|^{\alpha}.$$
From  $\Im \tilde{p}_\psi=0$ and \eqref{2.10}, we have
\begin{equation*}
\begin{array}{rl}
\varepsilon_1 |1+i\tau|^{\alpha}\leq& |\Im (1+i\tau)^\alpha|\leq 2|\sigma\tilde{X}\Sigma_{j=1}^{n}a_{jn}\xi_j|+O(\sqrt{X}(|\xi|+|\sigma|X)^2)\\
\leq& C|\sigma\tilde{X}|^2+\frac{\varepsilon_1}{4} |1+i\tau|^\alpha
\end{array}
\end{equation*}with some  constant $C>0$
which implies \eqref{2.11} and hence
\begin{equation}\label{2.12}
\begin{array}{l}
\{\Re \tilde{p}_\psi, \Im \tilde{p}_\psi\}_p\gtrsim (|1+i\tau|^\alpha+|\xi|^2+\sigma^2)^{3/2}.
\end{array}
\end{equation}

\section{Subelliptic estimates}\label{sec3}
In this section we will show the following subelliptic estimate for the operator $P_\psi$.
\begin{lemma}\label{lem3.1}
There exists a small constant $z_0$
such that for all $u(t,x,z)\in C_0^{\infty}(U\times[-z_0,z_0])\cap\dot{\mathcal{S}}(\bar{\mathbf{R}}_+^{1+n+1})$, we have that
\begin{equation}\label{3.1}
\begin{cases}
\begin{array}{l}
\Sigma_{k+s<2}||h(D_z)^{2-k-s}\Lambda_\alpha^{s}D_z^k u||+||h(D_z)\Lambda_\alpha u||\lesssim||P_\psi u||\quad {\rm if}\quad \alpha< \frac{4}{3},\\
\Sigma_{k+s<2}||h(D_z)^{2-k-s}\Lambda_\alpha^{s}D_z^k u||+||h(D_z)^{\frac{3}{2}-\frac{2}{\alpha}}\Lambda_\alpha^{\frac{2}{\alpha}} u||\lesssim||P_\psi u||\quad {\rm if}\quad \alpha\geq \frac{4}{3},
\end{array}
\end{cases}
\end{equation}
where $h(D_z)=(1+D_z^2)^{1/4}$, $U$ is an open neighborhood of the origin in ${\Bbb R}^{1+n}$ and we recall
$\Lambda_{\alpha}^s(\tau,\xi)=((1+|\xi|^2)^{1/\alpha}+i\tau)^{s\alpha/2}$.
\end{lemma}
\pf.
We will show
\begin{equation}\label{3.2}
\begin{array}{l}
|| u||^2_{H^{\frac{3\alpha}{4},\frac{3}{2}}}
\lesssim||P_\psi u||_{L^2}^2+||u||^2_{H^{\alpha/2,1}}.
\end{array}
\end{equation}
Here we have abused the notation $\Vert u\Vert_{H^{m,s}}$ with $m,\,s\in\mathbf{R}$ to denote the norm for
$u(t,x,z)\in\dot{\mathcal{S}}(\bar{\mathbf{R}}_+^{1+n+1})$ given by
$$
\Vert u\Vert_{H^{m,s}}^2:=\int\int\int(1+|\tau|^m+|\xi|^s+|\sigma|^s)^2|\hat u|^2\,d\tau\,d\xi\,d\sigma.
$$
By the fact given as Lemma 2.2 in \cite{Taylor}, p. 363 that for any small $\epsilon>0$,
there exists a small $z_0$ such that
\begin{equation}\label{3.3}
\begin{array}{l}
||u||^2_{H^{\frac{\alpha}{2},1}}
\lesssim \epsilon||h(D_z)u||^2_{H^{\frac{\alpha}{2},1}}.
\end{array}
\end{equation}
Hence, we have
\begin{equation}\label{3.4}
\begin{array}{l}
||u||^2_{H^{\frac{3\alpha}{4},\frac{3}{2}}}
\lesssim||P_\psi u||_{L^2}^2
\end{array}
\end{equation}
which implies \eqref{3.1} for a small enough $z_0$.

We write $||P_\psi u||_{L^2}^2=(P_\psi u,P_\psi u)$ as the following.
\begin{equation}\label{3.5}
\begin{array}{rl}
||P_\psi u||_{L^2}^2&=(P_\psi u,P_\psi u)\\
&=(P^*_\psi P_\psi u, u)\\
&=(P_\psi P^*_\psi u, u)+([P^*_\psi, P_\psi] u, u)\\
&=\big((I-\varpi B)P_\psi P^*_\psi u, u\big)+\big(([P^*_\psi, P_\psi]+\varpi BP_\psi P^*_\psi) u, u\big),
\end{array}
\end{equation}
where the symbol of the principal part of the commutator $[P^*_\psi, P_\psi]$
is $[\overline{\tilde{p}_\psi}, \tilde{p}_\psi]=2\{\Re \tilde{p}_\psi, \Im \tilde{p}_\psi\}$
which has been already studied in Section \ref{sec2}, $\varpi$ a large positive constant
and $B=\Lambda^{-1}$ with an elliptic pseudo-differential operator $\Lambda$ whose principal symbol is
$(|\tau|^\alpha+|\xi|^2+\sigma^2)^{1/2}$.
From \eqref{2.8}, we have that the principal symbol of $[P^*_\psi, P_\psi]+\varpi BP_\psi P^*_\psi$ satisfies
\begin{equation}\label{3.6}
\begin{array}{l}
\varpi(|\xi|^2+\sigma^2+|\tau|^\alpha)^{-1/2}|\tilde{p}_\psi|^2
+2\{\Re \tilde{p}_\psi, \Im \tilde{p}_\psi\}\gtrsim (|\xi|^2+\sigma^2+|\tau|^\alpha)^{3/2}.
\end{array}
\end{equation}
Then by G{\aa}rding's inequality and \eqref{3.5}, we obtain  \eqref{3.2} for large enough $\varpi$.

\section{Carleman estimates}\label{sec4}

In this section, we will derive a Carleman estimate from \eqref{3.1} by conjugating $u$ in \eqref{3.1} by $e^{i\beta z}$ with
the large parameter $\beta$.

\begin{lemma}\label{lem4.1}
There exist a sufficiently large constant $\beta_1>0$ depending on $n$
such that for all $v(t,x)\in C_0^{\infty}(U)\cap\dot{\mathcal{S}}(\bar{\mathbf{R}}_+^{1+n})$
and $\beta\geq \beta_1$, we have that
\begin{equation}\label{4.1}
\begin{cases}
\begin{array}{l}
\sum_{|\gamma|\leq1}\beta^{3-2|\gamma|}\int e^{2\beta\psi(x)}|D_x^\gamma v|^2dtdx
\lesssim \int e^{2\beta\psi(x)}|P(t,x,D_t,D_x) v|^2dtdx\quad \alpha<\frac{4}{3}\\
\sum_{|\gamma|\leq1}\beta^{3-2|\gamma|}\int e^{2\beta\psi(x)}|D_{x}^\gamma v|^2dtdx
+\beta^{3-\frac{4}{\alpha}}\int e^{2\beta\psi(x)}|D_{t}v|^2dtdx\\
\lesssim \int e^{2\beta\psi(x)}|P(t,x,D_t,D_x) v|^2dtdx\quad \alpha\geq\frac{4}{3},
\end{array}
\end{cases}
\end{equation}
where $\psi=\frac{1}{2}(x_n-2X)^2$.
\end{lemma}

\pf. Let $u\in C^\infty_0(U\times(-z_0,z_0))$ and $\beta>0$. We denote
$$u\hat{}(\sigma,t,x)=\int u(z,t,x)e^{-i\sigma z}dz.$$
Then
\begin{equation}\label{4.2}
\begin{array}{rl}
e^{-i\beta z}P_\psi (e^{i\beta z}u) =\frac{1}{2\pi}
\int e^{i(\sigma-\beta)z}p(t,x,D_t,D_x+i|\sigma|\nabla\psi)u\hat{}(\sigma-\beta,t,x)d\sigma.
\end{array}
\end{equation}
By the Leibniz formula, it yields that for any smooth function $w=w(t,x)$,
\begin{equation}\label{4.3}
\begin{array}{rl}
&p(t,x,D_t,D_x+i|\sigma|\nabla\psi)w\\
=&e^{|\sigma|\psi}p(t,x,D_t,D_x) e^{-|\sigma|\psi}w\\
=&e^{(|\sigma|-|\beta|)\psi}e^{|\beta|\psi}p(t,x,D_t,D_x)e^{-|\beta|\psi}e^{-(|\sigma|-|\beta|)\psi}w\\
=&p(t,x,D_t,D_x+i|\beta|\nabla\psi)w+\Sigma_{j+k+|\gamma|\le2,j>0}C_{j,k,\gamma}(t,x)(|\sigma|-|\beta|)^j|\beta|^kD_x^\gamma w\\
&+\Sigma_{l=1}^mC_l(|\sigma|-|\beta|)D_t^{\alpha_l-1} w,
\end{array}
\end{equation}
where the last term $\Sigma_{l=1}^mC_l(|\sigma|-|\beta|)D_t^{\alpha_l-1} w$ is coming from the last term of \eqref{2.1} whose symbol is $\Sigma_{l=1}^mq_l\frac{X}{T}i^\alpha(\tau-i)^{\alpha_l-1}\xi_n$.

Now, let $u(z,t,x)=f(t,x)g(z)$, then we have $u\hat{}(\sigma-\beta,t,x)=f(t,x)\hat{g}(\sigma-\beta)$.
Applying \eqref{4.3} to $w=f(t,x)\hat{g}(\sigma-\beta)$, we have from \eqref{4.2} and \eqref{4.3}
\begin{equation}\label{4.4}
\begin{array}{rl}
&e^{-i\beta z}P_\psi (e^{i\beta z}f(t,x)g(z))\\
=&\int e^{i(\sigma-\beta)z}p(t,x,D_t,D_x+i|\sigma|\nabla\psi)f(t,x)\hat{g}(\sigma-\beta)d\sigma\\
=&\int e^{i(\sigma-\beta)z}p(t,x,D_t,D_x+i|\beta|\nabla\psi)f(t,x)\hat{g}(\sigma-\beta)d\sigma\\
&+\int e^{i(\sigma-\beta)z}\Sigma_{j+k+|\gamma|\le 2,\,j>0}C_{j,k,\gamma}(t,x)(|\sigma|-|\beta|)^j|\beta|^kD_x^\gamma f(t,x)\hat{g}(\sigma-\beta)d\sigma\\
&+\Sigma_{l=1}^mC_l\int e^{i(\sigma-\beta)z}(|\sigma|-|\beta|)D_t^{\alpha_l-1} f(t,x)\hat{g}(\sigma-\beta)d\sigma\\
=&g(z)p(t,x,D_t,D_x+i|\beta|\nabla\psi)f(t,x)\\
&+\Sigma_{j+k+|\gamma|\le 2,\,j>0}C_{j,k,\gamma}(x)G_j(\beta)g(z)|\beta|^kD_x^\gamma f(t,x)\\
&+\Sigma_{l=1}^mC_lD_t^{\alpha_l-1} f(t,x)G_1(\beta)g(z),
\end{array}
\end{equation}
where $G_j(\beta)g(z)=\int e^{i(\sigma-\beta)z}(|\sigma|-|\beta|)^j\hat{g}(\sigma-\beta)d\sigma
=\int e^{i\sigma z}(|\sigma+\beta|-|\beta|)^j\hat{g}(\sigma)d\sigma$.
By the Plancherel theorem, we have
\begin{equation}\label{4.5}
\begin{array}{l}
||G_j(\beta)g(z)||^2\lesssim \int |\sigma|^{2j}|\hat{g}(\sigma)|^2d\sigma\lesssim||g||^2_{H^j(\mathbf{R})}.
\end{array}
\end{equation}
Let $g\in C_0^\infty((-z_0,z_0))$ be any fixed non-zero function. Combining \eqref{4.4} and \eqref{4.5}, we have
\begin{equation}\label{4.6}
\begin{array}{rl}
&||P_\psi (e^{i\beta z}f(t,x)g(z))||\\
\lesssim&||g||\cdot||p(t,x,D_t,D_x+i|\beta|\nabla\psi)f||+\Sigma_{j+k+|\gamma|\le 2,j>0}|\beta|^k\cdot||g||_{H^j(\mathbf{R})}\cdot||D_x^\gamma f||\\
&+\Sigma_{l=1}^m||D_t^{\alpha_l-1} f||\cdot||g||_{H^1(\mathbf{R})}\\
\lesssim&||p(t,x,D_t,D_x+i|\beta|\nabla\psi)f||+\Sigma_{k+|\gamma|<2}\beta^k\cdot||D_x^\gamma f||+||D_t^{\alpha-1} f||.
\end{array}
\end{equation}

On the other hand, we need to estimate the lower bound on the left hand side of \eqref{3.1}.
We let $u(t,x,z)=e^{i\beta z}f(t,x)g(z)$ in \eqref{3.1}.
A direct computation gives for $j\in \mathbf{N}$ and $k\in \mathbf{Z}_+:=\mathbf{N}\bigcup \{0\}$ that
\begin{equation}\label{4.7}
\begin{array}{rl}
&e^{-i\beta z}h(D_z)^jD_z^k(e^{i\beta z}g(z))\\
=&(2\pi)^{-1}\int e^{i(\sigma-\beta)z}h(\sigma)^j\sigma^k\hat{g}(\sigma-\beta)d\sigma\\
=&(2\pi)^{-1}\int e^{i\sigma z}h(\sigma+\beta)^j(\sigma+\beta)^k\hat{g}(\sigma)d\sigma\\
=&(2\pi)^{-1}\int e^{i\sigma z}h(\sigma+\beta)^j[\beta^k+\Sigma_{k'=0}^{k-1}C_{k'}\beta^{k'}\sigma^{k-k'}]\hat{g}(\sigma)d\sigma\\
=&(2\pi)^{-1}\int e^{i\sigma z}h(\sigma+\beta)^j\beta^k\hat{g}(\sigma)d\sigma
+(2\pi)^{-1}\Sigma_{k'=0}^{k-1}C_{k'}\int e^{i\sigma z}h(\sigma+\beta)^j\beta^{k'}\sigma^{k-k'}\hat{g}(\sigma)d\sigma\\
\end{array}
\end{equation}
and
\begin{equation}\label{4.8}
\begin{array}{rl}
&(2\pi)^{-1}\int e^{i\sigma z}h(\sigma+\beta)^j\beta^k\sigma^l\hat{g}(\sigma)d\sigma\\
=&\beta^k(2\pi)^{-1}\int e^{i\sigma z}[h(\sigma+\beta)-h(\beta)+h(\beta)]^j\sigma^l\hat{g}(\sigma)d\sigma\\
=&\beta^k(2\pi)^{-1}\int e^{i\sigma z}h(\beta)^j\sigma^l\hat{g}(\sigma)d\sigma+\beta^k H_{j,l}(\beta)g(z)
\end{array}
\end{equation}
with $H_{j,l}(\beta)g(z)
=(2\pi)^{-1}\int e^{i\sigma z}
\sum_{j'=0}^{j-1}(h(\sigma+\beta)-h(\beta))^{j-j'} h(\beta)^{j'}\sigma^{l}\hat{g}(\sigma)d\sigma$.

By \eqref{4.7} and \eqref{4.8}, we obtain
\begin{equation}\label{4.9}
\begin{array}{rl}
&e^{-i\beta z}h(D_z)^jD_z^k(e^{i\beta z}g(z))-h(\beta)^j\beta^kg(z)\\
=&\beta^kH_{j,0}(\beta)g(z)+
(2\pi)^{-1}\Sigma_{k'=0}^{k-1}C_{k'}\int e^{i\sigma z}h(\sigma+\beta)^j\beta^{k'}\sigma^{k-k'}\hat{g}(\sigma)d\sigma\\
=&\beta^kH_{j,0}(\beta)g(z)+\Sigma_{k'=0}^{k-1}C_{k'}\beta^{k'}H_{j,k-k'}(\beta)g(z)\\
&+(2\pi)^{-1}h(\beta)^j\Sigma_{k'=0}^{k-1}C_{k'}\beta^{k'}\int e^{i\sigma z}\sigma^{k-k'}\hat{g}(\sigma)d\sigma.
\end{array}
\end{equation}
From $|h(\sigma+\beta)-h(\beta)|\lesssim|\sigma|$, we can easily get
\begin{equation}\label{4.10}
\begin{array}{l}
||H_{j,l}(\beta)g(z)||\lesssim h(\beta)^{j-1}.
\end{array}
\end{equation}
Combining \eqref{4.9} and \eqref{4.10}, we have for large $\beta$ that
\begin{equation}\label{4.11}
\begin{array}{l}
\begin{cases}
||e^{-i\beta z}h(D_z)^jD_z^k(e^{i\beta z}g(z))-h(\beta)^j\beta^kg(z)||\lesssim h(\beta)^j\beta^{k}(h(\beta)^{-1}+\beta^{-1})\\
||h(D_z)^jD_z^k(e^{i\beta z}g(z))||\gtrsim h(\beta)^{j}\beta^{k},
\end{cases}
\end{array}
\end{equation}
By \eqref{4.9}, we have for large $\beta$
\begin{equation}\label{4.12}
\begin{array}{rl}
&||h(D_z)^j\Lambda_\alpha^{s}  D_z^k u||\\
=&||\Lambda_{\alpha}^{s}  f(t,x)e^{-i\beta z}h(D_z)^jD_z^k(e^{i\beta z}g(z))||\\
=&||\Lambda_{\alpha}^{s}  f(t,x)[e^{-i\beta z}h(D_z)^jD_z^k(e^{i\beta z}g(z))-h(\beta)^j\beta^kg(z)+h(\beta)^j\beta^kg(z)]||\\
\gtrsim&h(\beta)^{j}\beta^{k}||\Lambda_{\alpha}^{s}  f||.
\end{array}
\end{equation}
Recall that $h(\beta)\simeq\beta^{1/2}$, by
\eqref{4.12} and \eqref{3.1}, we obtain
\begin{equation}\label{4.13}
\begin{array}{rl}
&\sum_{|\gamma|\leq1}\beta^{3-2|\gamma|}\int |D_x^\gamma f|^2dtdx\\
\lesssim&\Sigma_{k+s<2}h(\beta)^{2(2-k-s)}\beta^{2k}||\Lambda_{\alpha}^{s} f||^2\\
\lesssim&\Sigma_{k+s<2}||h(D_z)^{2-k-s}\Lambda_{\alpha}^{s} D_z^k u ||^2\\
\lesssim&||P_\psi u||^2\\
=&||P_\psi(e^{i\beta z}f(t,x)g(z))||^2
\end{array}
\end{equation}
and
\begin{equation}\label{4.14}
\begin{cases}
\begin{array}{l}
\beta^{1/2}||D_t^{\alpha-1} f||\lesssim||P_\psi(e^{i\beta z}f(t,x)g(z))||,\quad \alpha<\frac{4}{3}\\
\beta^{1/2}||D_t^{\alpha-1} f||+\beta^{\frac{3}{2}-\frac{2}{\alpha}}||D_t f||\lesssim||P_\psi(e^{i\beta z}f(t,x)g(z))||,\quad \alpha\geq\frac{4}{3}.
\end{array}
\end{cases}
\end{equation}

Combining \eqref{4.13}, \eqref{4.14} and \eqref{4.6}, we have for large enough $\beta$ that
\begin{equation}\label{4.15}
\begin{cases}
\begin{array}{l}
\sum_{|\gamma|\leq1}\beta^{3-2|\gamma|}\int |D_x^\gamma f|^2dtdx\lesssim
||p(t,x,D_t,D_x+i\beta\nabla\psi)f||^2,\quad \alpha<\frac{4}{3}\\
\sum_{|\gamma|\leq1}\beta^{3-2|\gamma|}\int |D_x^\gamma f|^2dtdx+\beta^{3-\frac{4}{\alpha}}||D_t f||^2
\lesssim||p(t,x,D_t,D_x+i\beta\nabla\psi)f||^2,\quad \alpha\geq\frac{4}{3}.
\end{array}
\end{cases}
\end{equation}
By letting $f=e^{\beta \psi}v$ in \eqref{4.15}, we immediately have \eqref{4.1}.

\section{Proof of Theorem \ref{thm1.2}}\label{sec5}

This section is devoted to the proof of the Theorem \ref{thm1.2}.
The limit arguments in \eqref{4.1} imply that we can assume $u\in \dot{\mathcal{S}}(\overline{{\mathbf{R}}_+^{1+n}})$.
To apply Lemma \ref{lem4.1}, we first recall
$${\rm supp}\, u\subset\{x_n\geq 0\}$$
and then define a smooth function $\chi$ by

\begin{equation}\label{5.1}
\chi (x_n)=
\begin{cases}
\begin{array}{l}
1,\quad x_n\leq (1-\varepsilon)X,\\
0,\quad x_n\geq  X,
\end{array}
\end{cases}
\end{equation}
where $\varepsilon<1$ is a positive constant.
Recalling \eqref{1.7}, it is not hard to get
$(\chi u)(t,x)\in C_0^{\infty}(U)\cap\dot{\mathcal{S}}(\overline{{\mathbf{R}}_+^{1+n}})$.
Thus, we can apply the Carleman estimates \eqref{4.1} to $(\chi u)(t,x)$ and get that
\begin{equation}\label{5.2}
\begin{array}{rl}
&\sum_{|\gamma|\leq1}\beta^{3-2|\gamma|}\int_{x_n\leq (1-\varepsilon) X} e^{2\beta\psi(x)}|D^\gamma u|^2dtdx\\
\leq &\sum_{|\gamma|\leq1}\beta^{3-2|\gamma|}\int e^{2\beta\psi(x)}|D^\gamma (\chi u)|^2dtdx\\
\lesssim &\int e^{2\beta\psi(x)}|P(t,x,D_t,D_x) (\chi u)|^2dtdx\\
\lesssim &\sum_{|\gamma|\leq1}\int e^{2\beta\psi(x)}|\chi D^\gamma u|^2dtdx+
\int_{(1-\varepsilon) X<x_n\leq  X} e^{2\beta\psi(x)}|[P,\chi] u|^2dtdx,
\end{array}
\end{equation}
where $[\cdot,\cdot]$ denotes the commutator. Let $\beta$ be large enough to absorb the first term on the right hand side of \eqref{5.2} into the first and second lines, we get from \eqref{5.2} that
\begin{equation}\label{5.3}
\begin{array}{rl}
&\beta^{3}\int_{x_n\leq (1-2\varepsilon) X} e^{\beta (1+2\varepsilon)^2 X^2/2}|u|^2dtdx\\
\leq &\sum_{|\alpha|\leq1}\int_{(1-\varepsilon) X<x_n\leq  X} e^{2\beta\psi(x)}|D^\alpha u|^2dtdx\\
\lesssim &C(u)e^{\beta (1+\varepsilon)^2 X^2/2}.
\end{array}
\end{equation}
Let $\beta$ tend to $\infty$, we obtain that $u=0$ on $x_n\leq (1-2\varepsilon) X$.
Since this argument works for any $1>\varepsilon>0$, we get that
$u=0$ in $\{x_n< X\}$ which implies the conclusion of Theorem \ref{thm1.2}.

\section{Proof of Theorem \ref{thm1.1}}\label{sec6}
Without loss of generality, we may assume $\Omega=B(0,1)$ and $\omega=B(0,r)$ with $r<1$ and $u(t,y)=0$ in $(0,T)\times\omega$, where
$B(0,r)$ is a ball  in $\mathbf{R}^n$ with radius $r$ centered at the origin.
Let $y_r\in \partial\omega$, we can rotate
the coordinates such that the line passing through the origin and $y_r$ is on the $y_n$ axis.
We keep the same coordinates $y$ since it will not affect \eqref{1.3}.
For simplicity, we use a diffeomorphism which allows us to assume that $\Omega=(-1,1)^n$ and $\omega=(-r,r)^{n}$, respectively.
We keep the same coordinates $y$ since it only change the constant $\delta$ in \eqref{1.3}.

Now, we define a diffeomorphism $\tilde{y}$ from $\Omega$ to $\mathbf{R}^n$ by
$$\tilde{y}(y)=(\frac{y_1}{\sqrt{1-y_1^2}},\cdots,\frac{y_n}{\sqrt{1-y_n^2}}).$$
It is not hard to get that
$$\tilde{\emph{L}} u(t,\tilde{y}):=\emph{L} u(t,y)=\sum_{j,k=1}^n\tilde{a}_{j,k}(t,\tilde{y})\partial_{\tilde{y}_j}\partial_{\tilde{y}_k}u(t,\tilde{y})+{l}_1(t,\tilde{y};\nabla_y)u(t,\tilde{y}),$$
where $\tilde{a}_{j,k}(t,\tilde{y})=a_{j,k}(t, {y})(1+\tilde{y}_j^2)^{3/2}(1+\tilde{y}_k^2)^{3/2}$. From \eqref{1.3}, we have $\tilde{a}_{j,k}(t,\tilde{y})$ satisfying
\begin{equation}\label{6.1}
\begin{array}{l}
\delta |\tilde{\eta}|^2 \leq \sum_{j,k=1}^n\tilde{a}_{j,k}(t,\tilde{y})\eta_j\eta_k \leq \delta^{-1} |\tilde{\eta}|^2,
\end{array}
\end{equation}
where $\tilde{\eta}=((1+\tilde{y}_1^2)^{3/2}\eta_1,\cdots,(1+\tilde{y}_n^2)^{3/2}\eta_n)$
and $\eta_j$ corresponds to the symbol of $D_j:=-i\partial_{\tilde{y}_j}$.
Thus, by this diffeomorphism, $\Omega$ and $\omega$ are mapped to  $\Omega=\mathbf{R}^n$
and a rectangle $\tilde{\omega}=(-\tilde{r},\tilde{r})^{n}$ with some $\tilde{r}>0$, respectively.
To match the situation of Theorem \ref{thm1.2}, we note $\tilde{y}=\tilde{y}-y_{\tilde{r}}$,
where $y_{\tilde{r}}=(0,\cdots,0,\tilde{r})$.
\begin{equation}\label{6.2}
\left\{
\begin{array}{l}
\Sigma_{j=1}^mq_j\partial_t^{\alpha_j} u(t,\tilde y)-\tilde{\emph{L}} u(t,\tilde{y})=\tilde{\ell}_1(t,\tilde{y};\nabla_{\tilde y})u(t,\tilde{y}),\\
u(t,\tilde{y})=0,\,\,t\le0,\\
u(t,\tilde{y})=0,\,\,\tilde{y}\in\tilde{\omega}-y_{\tilde{r}},\,t\in(0,T).
\end{array}
\right.
\end{equation}

 It suffices to prove that
\begin{equation}\label{6.3}
\begin{array}{l}
u(t,\tilde{y})=0\quad{\rm in} \quad \cup_{s\in \mathbf{N}} E_s,
\end{array}
\end{equation}
where
$$E_s=\{(t,\tilde{y}):t\in(0,T),\, \tilde{y}_n+\frac{sX}{T}t<sX,\, |\tilde{y}'|<\sqrt{X}\}.$$
It should be noted that our arguments works for any $y_{\tilde{r}}\in \partial\tilde{\omega}$.

We will prove \eqref{6.3} by induction on $s$. Namely, we prove for $s=1$ and then we prove for any $s\ge 2$, $u=0$ in $E_{s-1}$ implies $u=0$ in $E_{s}$.
For $s=1$, we
apply the argument of the proof of Theorem \ref{thm1.2} for each $y_{\tilde{r}}\in \partial\tilde{\omega}$.
Namely, by using the change of variables
$x'=\tilde{y}', x_{1,n}:=x_n=\tilde{y}_n+c|\tilde{y}'|^2+\frac{X}{T}t, \tilde{t}=t$ with suitable constants
$c\geq 1$, $0<X\ll 1$, and using a cutoff function $\kappa_1$ such that $\kappa_1 u=u$ in $Q_1$, we will allow to have
\begin{equation}\label{6.4}
|x'|\lesssim \sqrt{X}, \quad x_n\leq X.
\end{equation}
Then the argument of proof of Theorem \ref{thm1.2} can give that $u=0$ in
$\{x_{1,n}=\tilde{y}_n+c|\tilde{y}'|^2+\frac{X}{T}t<X\}$.
Since the proof works for any $y_{\tilde{r}}$, this yields that  $u=0$ in $E_1$.

Next, we will prove for $s=2$ and each $y_{\tilde{r}}\in \partial\tilde{\omega}$.
Let us make the change of variables
$x'=\tilde{y}', x_{2,n}=\tilde{y}_n+c|\tilde{y}'|^2+\frac{2X}{T}t-X, \tilde{t}=t$.
It should be noted that $x_{2,n}=x_{1,n}+\frac{X}{T}t-X$.
Recall
\begin{equation*}
\begin{array}{l}
\begin{cases}
Q_\zeta=\{(\tilde{y}',\tilde{y}_n):|\tilde{y}'|\leq \sqrt{X},\,\,-l/3<\tilde{y}_n\leq \zeta X \},\\
\tilde{Q}_\zeta=\{(\tilde{y}',\tilde{y}_n):|\tilde{y}'|\leq 2\sqrt{X},\,\,-2l/3<\tilde{y}_n\leq (\zeta+1) X \}
\end{cases}
\end{array}
\end{equation*}
and
\begin{equation*}
\begin{array}{l}
\kappa_\zeta(\tilde{y})=
\begin{cases}
1, \quad \tilde{y}\in Q_\zeta,\\
0, \quad \tilde{y}\in \mathbf{R}^n\setminus \tilde{Q}_\zeta.
\end{cases}
\end{array}
\end{equation*}

From  $u=0$ in $\{x_{1,n}=\tilde{y}_n+\frac{X}{T}t<X\}$ for any $y_{\tilde{r}}\in \partial\tilde{\omega}$
and the definition of $\kappa_2(y)$, we have
$${\rm supp}\, \kappa_2 u\subset\{x_{2,n}\geq 0\}.$$
For this, if we assume that $x_{2,n}<0$, then $x_{1,n}<X$ which implies $u=0$ in $x_{2,n}<0$.
Since \begin{equation*}
\chi (x_{2,n})=
\begin{cases}
\begin{array}{l}
1,\quad x_{2,n}\leq (1-\epsilon)X,\\
0,\quad x_{2,n}\geq X,
\end{array}
\end{cases}
\end{equation*}
we get that
$${\rm supp}\, \chi\kappa_2 u\subset Q_2.$$
To see this, we can assume that $\tilde{y}_n>2X$ then $x_{2,n}>X$ which implies $\chi=0$ in $\tilde{y}_n>2X$.
If we assume that $|\tilde{y}'|> \sqrt{X}$
then $x_{1,n}<X$ which implies $u=0$ in $x_{1,n}<X$.
It is easy to see that $\kappa_2 u=u$ in $Q_2$. Thus, we can denote $\kappa_2 u$ by $u$ and assume that
\begin{equation}\label{6.5}
|x'|\lesssim \sqrt{X}, \quad x_{2,n}\leq X.
\end{equation}
Repeat the arguments in Theorem \ref{thm1.2} for each $y_{\tilde{r}}\in \partial\tilde{\omega}$,
we can get that $u=0$ in $E_2$.

By induction, we assume that  $u=0$ in $\{x_{s-1,n}=\tilde{y}_n+\frac{(s-1)X}{T}t-(s-2)X<X\}$.
We use the change of variables $x'=\tilde{y}', x_{s,n}=\tilde{y}_n+c|\tilde{y}'|^2+\frac{sX}{T}t-(s-1)X, \tilde{t}=t$.
By considering the support of $\chi\kappa_s u\subset Q_s$, we can assume as for the previous arguments that
\begin{equation}\label{6.6}
|x'|\lesssim \sqrt{X}, \quad x_{s,n}\leq X.
\end{equation}
Here note that $|\tilde{y}_n|\le s X$. To apply the previous arguments,
we need to calarify how the estimate of Poisson bracket depends on $s$. From \eqref{6.2}, the principal operator is
\begin{equation*}
\begin{array}{l}
\partial_t^\alpha-\sum_{j,k=1}^n\tilde{a}_{j,k}(t,\tilde{y})\partial_{\tilde{y}_j}\partial_{\tilde{y}_k}u
\end{array}
\end{equation*}
with $\tilde{a}_{j,k}(t,\tilde{y})=a_{j,k}(t, {y})(1+\tilde{y}_j^2)^{3/2}(1+\tilde{y}_k^2)^{3/2}$.

In the new coordinates $(t,\tilde{y})$, we use the Holmgren type transform
$x'=\tilde{y}', x_{s,n}=\tilde{y}_n+c|\tilde{y}'|^2+\frac{sX}{T}t-(s-1)X, \tilde{t}=t$. By abusing the notation $P_\psi$, we denote by $P_\psi$ the corresponding one for the operator on the left side of the equation in \eqref{6.2}.
We note that the principal symbol $\tilde{p}_\psi$ of $P_\psi$
is given by
\begin{equation*}
\begin{array}{rl}
\tilde{p}_\psi=&(1+i\tau)^{\alpha}+\tilde{a}_{nn}(\eta_n+i|\sigma|\hat{X})^2\\
&+2\Sigma_{j=1}^{n-1}\tilde{a}_{jn}(\eta_j+2cx_j\eta_n+2icx_j|\sigma|\hat{X})(\eta_n+i|\sigma|\hat{X})\\
&+\Sigma_{j,k=1}^{n-1}\tilde{a}_{jk}(\eta_j+2cx_j\xi_n+2icx_j|\sigma|\hat{X})(\eta_k+2cx_k\eta_n+2icx_k|\sigma|\hat{X}),
\end{array}
\end{equation*}
where $\psi=\frac{1}{2}(x_{s,n}-2X)^2$, $\hat{X}=x_{s,n}-2X$ and we denote $x_{s,n}$ by $x_n$,
and the symbol of each $D_{x_j}$ is $\eta_j$.

Recall the principal part $\{\Re \hat{p}, \Im \hat{p}\}_p$ of
the Poisson bracket $\{\Re \hat{p}, \Im \hat{p}\}$
\begin{equation*}
\begin{array}{ll}
\{\Re \tilde{p}_\psi, \Im \tilde{p}_\psi\}_p
=\Sigma_{j=1}^{n}(\partial_{\eta_j}\Re \tilde{p}_\psi\cdot\partial_{x_j}\Im \tilde{p}_\psi
-\partial_{x_j}\Re \tilde{p}_\psi\cdot\partial_{\eta_j}\Im \tilde{p}_\psi).
\end{array}
\end{equation*}
Write $\tilde{p}_\psi=\tilde{p}_{1,\psi}+\tilde{p}_{2,\psi}$, where
\begin{equation*}
\begin{array}{rl}
\tilde{p}_{1,\psi}=&\tilde{a}_{nn}(\eta^2_n-|\sigma|^2\hat{X}^2)+2\Sigma_{j=1}^{n-1}\tilde{a}_{jn}\eta_j\eta_n+\Sigma_{j,k=1}^{n-1}\tilde{a}_{jk}\eta_j\eta_k\\
&+i[2\tilde{a}_{nn}\eta_n|\sigma|\hat{X}+2\Sigma_{j=1}^{n-1}\tilde{a}_{jn}\eta_j|\sigma|\hat{X}]
\end{array}
\end{equation*}
and
\begin{equation*}
\begin{array}{rl}
\tilde{p}_{2,\psi}=&(1+i\tau)^\alpha+2\Sigma_{j=1}^{n-1}\tilde{a}_{jn}(2cx_j\eta_n+2icx_j|\sigma|\hat{X})(\eta_n+i|\sigma|\hat{X})\\
&+\Sigma_{j,k=1}^{n-1}\tilde{a}_{jk}(2cx_j\eta_n+2icx_j|\sigma|\tilde{X})(\eta_k+2cx_k\eta_n+2icx_k|\sigma|\hat{X})\\
&+2\Sigma_{j,k=1}^{n-1}\tilde{a}_{jk}\eta_j(2cx_k\eta_n+2icx_k|\sigma|\hat{X}).
\end{array}
\end{equation*}
A direct computation gives
\begin{equation*}
\begin{array}{rl}
\partial_{\eta_n}\Re \tilde{p}_{1,\psi}=&2\Sigma_{j=1}^{n}\tilde{a}_{jn}\eta_j\\
\partial_{x_n}\Im \tilde{p}_{1,\psi}=&2\partial_{x_n}(\Sigma_{j=1}^{n}\tilde{a}_{jn})\eta_j|\sigma|\hat{X}
+2|\sigma|\Sigma_{j=1}^{n}\tilde{a}_{jn}\eta_j\\
=&2|\sigma|\Sigma_{j=1}^{n}\tilde{a}_{jn}\eta_j+O((|\tilde{\eta}|+|\tilde{\sigma}X|)^2)\\
\partial_{x_n}\Re \tilde{p}_{1,\psi}=&(\eta^2_n-|\sigma|^2\hat{X}^2)\partial_{x_n}(\tilde{a}_{nn})+2\Sigma_{j=1}^{n-1}\partial_{x_n}(\tilde{a}_{jn})\eta_j\eta_n\\
&+\Sigma_{j,k=1}^{n-1}\partial_{x_n}(\tilde{a}_{jn})\eta_j\eta_k-2\tilde{a}_{nn}\sigma^2\hat{X}\\
=&-2\tilde{a}_{nn}\sigma^2\hat{X}+O((|\tilde{\eta}|+|\tilde{\sigma}X|)^2)\\
\partial_{\eta_n}\Im \tilde{p}_{1,\psi}=&2\tilde{a}_{nn}|\sigma|\hat{X}
\end{array}
\end{equation*}
and
\begin{equation*}
\begin{array}{rl}
\nabla_{\eta'}\Re \tilde{p}_{1,\psi}=&O(|\tilde{\eta}|)\\
\nabla_{x'}\Im \tilde{p}_{1,\psi}=&O((|\tilde{\eta}|+|\tilde{\sigma}X|)^2)\\
\nabla_{x'}\Re \tilde{p}_{1,\psi}=&O((|\tilde{\eta}|+|\tilde{\sigma}X|)^2)\\
\nabla_{\eta'}\Im \tilde{p}_{1,\psi}=&O(|\tilde{\sigma}X|)
\end{array}
\end{equation*}
where $\tilde{\sigma}=(1+\tilde{y}_n^2)^{3/2}\sigma$.

Thus, we have
\begin{equation}\label{6.7}
\begin{array}{rl}
\partial_{\eta_n}\Re \tilde{p}_{\psi}=&2\Sigma_{j=1}^{n}\tilde{a}_{jn}\eta_j+O(\sqrt{X}(|\tilde{\eta}|+|\tilde{\sigma}X|))\\
\partial_{x_n}\Im \tilde{p}_{\psi}=&2|\sigma|\Sigma_{j=1}^{n}\tilde{a}_{jn}\eta_j+O(\frac{1}{\sqrt{X}}(|\tilde{\eta}|+|\tilde{\sigma}X|)^2)\\
\partial_{x_n}\Re \tilde{p}_{\psi}=&-2\tilde{a}_{nn}\sigma^2\hat{X}+O(\frac{1}{\sqrt{X}}(|\tilde{\eta}|+|\tilde{\sigma}X|)^2)\\
\partial_{\eta_n}\Im \tilde{p}_{\psi}=&2\tilde{a}_{nn}|\sigma|\hat{X}+O(\sqrt{X}(|\tilde{\eta}|+|\tilde{\sigma}X|))
\end{array}
\end{equation}
and
\begin{equation}\label{6.8}
\begin{array}{rl}
\nabla_{\eta'}\Re \tilde{p}_{\psi}=&O(|\tilde{\eta}|+|\tilde{\sigma}X|)\\
\nabla_{x'}\Im \tilde{p}_{\psi}=&O((|\tilde{\eta}|+|\tilde{\sigma}X|)^2)\\
\nabla_{x'}\Re \tilde{p}_{\psi}=&O((|\tilde{\eta}|+|\tilde{\sigma}X|)^2)\\
\nabla_{\eta'}\Im \tilde{p}_{\psi}=&O(|\tilde{\eta}|+|\tilde{\sigma}X|).
\end{array}
\end{equation}
Combining \eqref{6.7} and \eqref{6.8}, we have
\begin{equation}\label{6.9}
\begin{array}{rl}
\{\Re \tilde{p}_\psi, \Im \tilde{p}_\psi\}_p=&\Sigma_{j=1}^{n}(\partial_{\eta_j}\Re \tilde{p}\cdot\partial_{x_j}\Im \tilde{p}
-\partial_{x_j}\Re \tilde{p}\cdot\partial_{\eta_j}\Im \tilde{p})\\
=&4|\sigma|(\Sigma_{j=1}^{n}\tilde{a}_{jn}\eta_j)^2+4\tilde{a}_{nn}^2|\sigma|^3\hat{X}^2+O(\frac{1}{\sqrt{X}}(|\tilde{\eta}|+|\tilde{\sigma}X|)^3)\\
\gtrsim &4|\sigma|(\Sigma_{j=1}^{n}\tilde{a}_{jn}\eta_j)^2+(1+\tilde{y}_n^2)^{3/2}X^{-1}|\tilde{\sigma}\hat{X}|^3
+O(\frac{1}{\sqrt{X}}(|\tilde{\eta}|+|\tilde{\sigma}X|)^3).
\end{array}
\end{equation}

The following Lemma needs the help of \eqref{6.1}
\begin{equation*}
\begin{array}{l}
\delta |\tilde{\eta}|^2 \leq \sum_{j,k=1}^n\tilde{a}_{j,k}(t,\tilde{y})\eta_j\eta_k \leq \delta^{-1} |\tilde{\eta}|^2,
\end{array}
\end{equation*}
where $\tilde{\eta}=((1+\tilde{y}_1^2)^{3/2}\eta_1,\cdots,(1+\tilde{y}_n^2)^{3/2}\eta_n)$.

\begin{lemma}\label{lem6.1}
Let $\tilde{p}_{\psi}=0$, we have that
\begin{equation}\label{6.10}
\begin{array}{l}
\{\Re \tilde{p}_\psi, \Im \tilde{p}_\psi\}_p
\gtrsim (1+\tilde{y}_n^2)^{3/2}(|\tilde{\eta}|^2+\tilde{\sigma}^2+|\tau|^\alpha)^{3/2}.
\end{array}
\end{equation}
\end{lemma}
\pf. Recall that
\begin{equation*}
\begin{array}{rl}
\tilde{p}_\psi=&\Sigma_{l=1}^mq_l(1+i\tau)^{\alpha_l}+\tilde{a}_{nn}(\eta_n+i|\sigma|\tilde{X})^2\\
&+2\Sigma_{j=1}^{n-1}\tilde{a}_{jn}(\eta_j+2x_j\eta_n+2ix_j|\sigma|\hat{X})(\eta_n+i|\sigma|\hat{X})\\
&+\Sigma_{j,k=1}^{n-1}\tilde{a}_{jk}(\eta_j+2x_j\eta_n+2ix_j|\sigma|\hat{X})(\eta_k+2x_k\eta_n+2ix_k|\sigma|\hat{X})\\
=&(1+i\tau)^\alpha+\Sigma_{j,k=1}^{n}\tilde{a}_{jk}\eta_j \eta_k - \tilde{a}_{nn}|\sigma|^2\hat{X}^2\\
&+2i|\sigma|\hat{X}\Sigma_{j=1}^{n}\tilde{a}_{jn}\eta_j+O(\sqrt{X}(|\tilde{\eta}|+|\tilde{\sigma}|X)^2).
\end{array}
\end{equation*}

Let us divide $\tau$ into two cases. If $|\tau| \leq 1$, then
\begin{equation}\label{6.11}
\begin{array}{l}
\Re (1+i\tau)^\alpha \geq \varepsilon_0 |1+i\tau|^\alpha
\end{array}
\end{equation}
where $\varepsilon_0$ is a positive constant.
From $\Re \tilde{p}_\psi=0$ and \eqref{6.11}, we have
\begin{equation*}
\begin{array}{l}
\Re\big(\Sigma_{l=1}^mq_l(1+i\tau)^{\alpha_l}\big)+\Sigma_{j,k=1}^{n}\tilde{a}_{jk}\eta_j \eta_k=\tilde{a}_{nn}|\sigma|^2\tilde{X}^2+O(\sqrt{X}(|\tilde{\eta}|+|\tilde{\sigma}|X)^2)
\end{array}
\end{equation*}
which implies
\begin{equation*}
\begin{array}{l}
|\tilde{\eta}|^2+|(1+i\tau)|^\alpha\lesssim |\tilde{\sigma}|^{2}X^2.
\end{array}
\end{equation*}

If  $|\tau| \geq 1$, then
\begin{equation}\label{6.12}
\begin{array}{l}
|\Im (1+i\tau)^\alpha| \geq \varepsilon_1 |1+i\tau|^\alpha
\end{array}
\end{equation}
where $\varepsilon_1$ is a positive constant.

First, we let
$$(2\Sigma_{l=1}^mq_l)|1+i\tau|^\alpha \leq \delta_0|\tilde{\eta}|^2,$$
where $\delta_0$ is the positive constant in \eqref{6.1} satisfying $\delta_0|\tilde{\eta}|^2\leq\Sigma_{j,k=1}^{n}\tilde{a}_{jk}\eta_j \eta_k$.
From  $\Re \tilde{p}_\psi=0$, we have
\begin{equation*}
\begin{array}{l}
\Re\big(\Sigma_{l=1}^mq_l(1+i\tau)^{\alpha_l}\big)+\Sigma_{j,k=1}^{n}\tilde{a}_{jk}\eta_j \eta_k=\tilde{a}_{nn}|\sigma|^2\tilde{X}^2+O(\sqrt{X}(|\tilde{\eta}|+|\tilde{\sigma}|X)^2)
\end{array}
\end{equation*}
which implies
\begin{equation*}
\begin{array}{l}
|\tilde{\eta}|^2+|(1+i\tau)|^\alpha\lesssim |\tilde{\sigma}|^{2}X^2.
\end{array}
\end{equation*}

On the other hand, we let
$$\delta_0|\tilde{\eta}|^2 \leq (2\Sigma_{l=1}^mq_l)|1+i\tau|^\alpha.$$
\begin{equation*}
\begin{array}{rl}
\varepsilon_1 |1+i\tau|^\alpha\leq& |\Im (1+i\tau)^\alpha|\leq 2|\sigma\tilde{X}\Sigma_{j=1}^{n}\tilde{a}_{jn}\eta_j|+O(\sqrt{X}(|\tilde{\eta}|+|\tilde{\sigma}|X)^2)\\
\leq& C|\sigma\tilde{X}|^2+\frac{\varepsilon_1}{4} |1+i\tau|^\alpha
\end{array}
\end{equation*}
which implies
\begin{equation*}
\begin{array}{l}
|\tilde{\eta}|^2+|(1+i\tau)|^\alpha\lesssim |\tilde{\sigma}|^{2}X^2.
\end{array}
\end{equation*}
Thus, we have finished proving \eqref{6.10} in Lemma \ref{lem6.1}.

We use \eqref{6.10}
to the proof of subelliptic estimate together with \eqref{3.5} and \eqref{3.6}.
Then the previous estimate \eqref{3.2}
changes to
\begin{equation}\label{6.13}
\begin{array}{l}
|| u||^2_{H^{\frac{3\alpha}{4}\frac{3}{2}}}
\lesssim||P_\psi u||_{L^2}^2+C(s)||u||^2_{H^{\alpha/2,1}}
\end{array}
\end{equation}
with some constant $C(s)>0$. Also the previous estimate \eqref{3.3}
changes to the following. That is there exists a small $z_0(s)>0$ depending on $s$ such that
\begin{equation}\label{6.14}
\begin{array}{l}
||u||^2_{H^{\frac{\alpha}{2},1}}
\lesssim \epsilon(s)||h(D_z)u||^2_{H^{\frac{\alpha}{2},1}}
\end{array}
\end{equation}
for $u=u(t,x,z)\in C_0^\infty(U\times[-z_0(s),z_0(s)])\cap \dot{\mathcal{S}}(\overline{\mathbb R}_+^{1+n+1})$,
where $\epsilon(s)$ is a small positive constant depending on $s$.

Then in the proof of the Carleman estimate \eqref{4.1},
the constants in the estimates from \eqref{4.6} to \eqref{4.15} depends on $s$
because the derivatives of $g\in C_0^\infty((-z_0(s), z_0(s))$ are involved in these estimates.
Therefore the Carleman estimate \eqref{4.1} will hold for any $\beta\ge\beta_1(s)$ with some $\beta_1(s)>0$.

Now by repeating the arguments in Theorem \ref{thm1.2} for each $y_{\tilde{r}}\in \partial\tilde{\omega}$,
we can get that $u=0$ in $E_s$, which completes our proof.

\end{document}